\documentclass[11pt]{amsart}
\usepackage{geometry}                % See geometry.pdf to learn the layout options. There are lots.
\usepackage{csquotes}
\usepackage{graphicx}
\usepackage{amssymb}
\usepackage{amsthm}
\usepackage{epstopdf}
\usepackage{color}
\usepackage{upgreek}
\label{csquotes}

\usepackage{url}
\usepackage[square,numbers]{natbib}
\usepackage{subeqnarray}
\definecolor{darkred}{rgb}{.6,0,0}
\definecolor{darkblue}{rgb}{0,0,.7}
\newtheorem{example}{Example}[section]
\newtheorem{remark}{Remark}[section]

\title[Online Digital Twins]{Online Digital Twins: Data Assimilation and Control}
\author{Sebastian Reich}
\address{Institut f\"ur Mathematik, Universit\"at Potsdam}
\begin{document}

\begin{abstract}
    The online assimilation of data and the application of controls are essential ingredients of the digital twin paradigm. In this paper, we formalize the online interaction between the physical process and its virtual digital representation in the setting of continuously arriving data from the physical process and a receding horizon virtual planning process, which provides the control for the physical process. Our mathematical formulation relies on a mean field representation of the data assimilation process as well as a mean field Pontryagin minimum principle for the optimal control problem. Such mean field formulations are computationally attractive because they allow for simultaneous implementation in terms of Hamiltonian interacting particle systems. A more detailed discussion is included for the ensemble Kalman--Bucy filter in combination with affine approximations of the receding horizon control laws. This paper sets out the mathematical framework, while a computational study is left for future research. 
\end{abstract}

\maketitle

%%%%%%%%%%%%%%%%%%%%%%%%%%%%%%
%
%%%%%%%%%%%%%%%%%%%%%%%%%%%%%%

\section{Introduction}

The two way interactions between a physical process and its virtual digital representation form the essential component of the digital twin paradigm, as discussed in detail in \cite{NA24}. Partial and noisy data are collected from the physical process and are used to continuously adjust the virtual representation through what is called filtering 
\cite{jazwinski2007stochastic,bain2008fundamentals} or data assimilation \cite{reich2015probabilistic,Evensenetal2022}. A data tuned virtual representation in the form of computational models allows in turn for online planning and the design of control laws which are applicable to the physical process. This planning step can be formalized as an optimal control problem for partially observed diffusion processes \cite{ASTROM1965,Bensoussan92,Bensoussan,Handel07}; these are also often referred to as partially observed Markov decision processes (POMDP) in the engineering literature \cite{Meyn}.

The connection of the digital twin paradigm to POMDP and probabilistic graphic models has been made before \cite{PGM,Willcox24}. This paper further develops this perspective in the context of continuous time data assimilation and control problems.  Furthermore, the approach followed in this paper is rooted in mean field formulations of both data assimilation and optimal control. On the data assimilation side, such formulations are well established and, for example, provided by the feedback particle filter (FPF) \cite{meyn13,reich2019data} and its approximation in the form of the ensemble Kalman--Bucy filter (EnKFB) \cite{bergemann2012ensemble,TdWMR17}. On the optimal control side, we rely on the recently proposed mean field Pontryagin minimum principle (PMP) \cite{R26a}. The proposed mean field formulation is Hamiltonian and facilitates a combined online treatment of data assimilation and control and, thus, avoids the separation of both tasks as inherent in most current digital twin implementations. See, for example, \cite{PGM}.

First attempts to apply mean field data assimilation algorithms to the control of partially observed diffusion processes can be found in the pioneering work \cite{PM13} in the context of infinite horizon optimal control formulations of partially observed diffusion processes. The formulation proposed in \cite{PM13} replaces infinite dimensional belief states by a finite but high dimensional state space; thus leading to a high dimensional Hamilton--Jacobi--Bellman equation in an approximate value function, which is still computationally intractable in general. More recently, \cite{R25c} combined the EnKBF and mean field PMP to provide a computationally tractable approximation which however ignores the impact of future observations. These contributions assume continuous time observations; a PMP framework for optimal control of finite horizon diffusion processes subject to discrete time observations has recently been put forward in \cite{BNST25}. Here, we focus on receding horizon control formulations; also called model predictive control (MPC) \cite{MPC}, and continuous time observations instead. An EnKBF based receding horizon optimal planning approach has already been considered in \cite{R25a} using the stochastic PMP principle leading to forward and backward stochastic differential equations \cite{Carmona}. The novel contribution of this paper is a Hamiltonian mean field and interacting particle formulation of receding horizon optimal planning under uncertainty in the context of digital twins. Numerical implementations will be considered in future work. 

The remainder of this paper is organized as follows. Section \ref{sec:problem} provides a mathematical specification of the digital twin paradigm in continuous time. Mean field data assimilation formulations and their interacting particle implementations are summarized in Section \ref{sec:DA}, while the mean field PMP formulation of \cite{R26a} is extended to state dependent diffusion matrices in Section \ref{sec:SOC}. Such an extension is crucial for allowing an extension of the FPF to the planning process in a receding horizon optimal control context. The general framework is provided in Section \ref{sec:POMDP} while Section \ref{sec:Example} contains a more detailed discussion in the context of the EnKFB and affine approximations to the control laws. The paper concludes with some remarks on computational challenges and possible extensions in Section \ref{sec:conclusions}.

%%%%%%%%%%%%%%%%%%%%%%%%%%%%%%
%
%%%%%%%%%%%%%%%%%%%%%%%%%%%%%%

\section{Problem formulation} \label{sec:problem}

In order to formalize the digital twin paradigm we adapt a continuous time setting and assume that the temporal physical process ${\rm Z}_t^\dagger$ with $t\ge 0$ denoting time can be steered through a control $U_t \in \mathbb{R}^{d_u}$. Control actions are based on partial observations $Y_t^\dagger \in \mathbb{R}^{d_y}$ of the physical process. We formalize this setting in the following manner. We assume that ${\rm Z}_t^\dagger$, $t>0$, takes values in a Banach space $\mathcal{P}$ and is unique to the initial ${\rm Z}_0^\dagger$ and the applied control sequence $U_{0:t} = \{U_s\}_{0\le s<t}$. We denote this dependency by ${\rm Z}_t^\dagger|{\rm Z}_0^\dagger,U_{0:t}$. Observations $Y_t^\dagger$ are assumed to satisfy the following differential equation:
\begin{equation} \label{eq:physical}
	\dot{Y}_t^\dagger = {\rm h}_t^\dagger({\rm Z}_t^\dagger|{\rm Z}_0^\dagger,U_{0:t}).
\end{equation}
Here, ${\rm h}_t^\dagger:\mathcal{P} \to \mathbb{R}^{d_y}$ is an appropriate, possibly time dependent function. Except for the observed $Y_{0:t}^\dagger = \{Y_s^\dagger\}_{0\le s < t}$ the physical process is inaccessible to the virtual environment. The control sequence $U_{0:t} = \{U_s\}_{0\le s < t}$, on the other hand, is provided by the virtual environment and is typically chosen so that a cost is minimized or a reward maximized.

The description (\ref{eq:physical}) can be seen in the tradition of what is called {\it Laplace's demon} of which Laplace himself wrote in \cite{Laplace}: 
\begin{displayquote}All these efforts in the search for truth tend to lead back continually to the vast intelligence which we have just mentioned, but from which it will always remain infinitely removed.\end{displayquote} 
Therefore, any virtual representation can only approximate (\ref{eq:physical}) and, for the purpose of this paper, we assume that the evolution model that describes a virtual representation of (\ref{eq:physical}) is of the form
\begin{subequations} \label{eq:virtual}
	\begin{align}
		\dot{Z}_t &= f(Z_t,U_t) + \Gamma \dot{V}_t,\\
		\dot{Y}_t &= h(X_t) + \dot{W}_t.
	\end{align}
\end{subequations}
Brownian motions $V_t \in \mathbb{R}^{d_z}$ and $W_t \in \mathbb{R}^{d_y}$ represent model and measurement errors, respectively, which mimic the discrepancy between the physical process and its virtual representation. Throughout this paper, we formally differentiate with respect to time, which should be interpreted in the sense of It\^o \cite{Pavliotis2016} unless otherwise stated.

The initial distribution $\rho_0$ of $Z_0$, the functions $f:\mathbb{R}^{d_z}\times \mathbb{R}^{d_u} \to \mathbb{R}^{d_z}$ and $h:\mathbb{R}^{d_z}\to \mathbb{R}^{d_y}$, and the matrix $\Gamma\in \mathbb{R}^{d_z\times d_z}$ completely characterize the temporal evolution of the virtual model for $t>0$ provided the control sequence $U_{0:t}$ is also given. We denote the induced law of $Z_t$ by $\rho_t$. We emphasize that the dimensions of ${\rm Z}_t^\dagger \in \mathcal{P}$ and $Z_t\in \mathbb{R}^{d_z}$ are different and $\dim (\mathcal{P})\gg d_z$, in general.

Assuming that a control sequence $U_{0:t_\ast}$ has been specified up to the present time $t_\ast>0$, the assimilation of the data $Y^\dagger_{0:t_\ast}$ leads to a change in measure from $Z_t \sim \rho_t$ to  $Z_t|Y^\dagger_{0:t} \sim \mu_t$ with the time evolution of the filtering distribution $\mu_t$ characterized by the Kushner--Stratonovich equation associated with (\ref{eq:virtual}) \cite{bain2008fundamentals}. We note that the derivation of the Kushner--Stratonovich equation assumes that the virtual representation (\ref{eq:virtual}) is a replica of the physical process (\ref{eq:physical}); {\it i.e.}, $\mathcal{P} = \mathbb{R}^{d_z}$ and
\begin{subequations}
	\begin{align}
		\dot{\rm Z}_t^\dagger &= f({\rm Z}_t^\dagger,U_t) + \Gamma
        \dot{V}_t^\dagger,\\
		\dot{Y}_t^\dagger &= h({\rm Z}_t^\dagger) + \dot{W}^\dagger_t
	\end{align}
\end{subequations}
with $W_t^\dagger \in \mathbb{R}^{d_y}$ and $V_t^\dagger \in \mathbb{R}^{d_z}$ denoting fixed realizations of Brownian motion. This setting is called the ideal twin scenario and amounts to the existence of {\it Laplace's demon}. Therefore, important challenges for the digital twin paradigm follow from discrepancies between physical process and virtual representation; {\it i.e.}, between (\ref{eq:physical}) and (\ref{eq:virtual}).

We have assumed that the control $U_t$ is known up to the present time $t_\ast$. Finding the optimal control beyond $t_\ast$ requires planning for the future $\tau > t_\ast$ with regard to a given cost function. This planning step must be initiated at a time $t = t_p < t_\ast$ and completed by time $t=t_\ast$. We formulate the planning step in the form of a receding horizon stochastic optimal control problem with cost function
\begin{equation} \label{eq:cost}
	J(U_{t_p:t_p+T};\upmu_{t_p}) = \mathbb{E}^\ddagger \left[ \int_{t_p}^{t_p+T}
	\left( \upmu_\tau[r] + \frac{1}{2} \|U_\tau\|_R^2\right) {\rm d}\tau + \upmu_{t_p+T}[g]
	\right],
\end{equation}
where the time interval $T>0$ is chosen such that $t_\ast \in (t_p,t_p+T)$. We assume a quadratic cost on the control for simplicity. Furthermore, $R \in \mathbb{R}^{d_u \times d_u}$ is a symmetric positive definite matrix and $\|u\|_{R}^2 = u^{\rm T}R^{-1}u$, $r(z)$ denotes the running cost and $g(z)$ the terminal cost,
\begin{equation}
	\upmu_\tau[r] = \int r(z) \upmu_\tau({\rm d}z), \quad
	\upmu_{t+T}[g] = \int g(z) \upmu_{t+T}({\rm d}z)
\end{equation}
their expectation values with respect to the law $\upmu_\tau$ of $Z_\tau$ conditioned on virtual observations, which we denote by $Y_{t_p:t_p+T}^\ddagger$. The expectation $\mathbb{E}^\ddagger$ in (\ref{eq:cost}) is taken with respect to these virtual observations.  

The conditional distributions $\upmu_\tau$ still satisfy the Kushner--Stratonovich equation with initial condition $\upmu_{t_p} = \mu_{t_p}$; but now driven by virtual observations $Y_{t_p:t_p+T}^\ddagger$ instead of an actual observation $Y_{t_p:t_p+T}^\dagger$ of the physical process. Note that $\upmu_\tau$ in (\ref{eq:cost}) are therefore random probability measures, which depend on the realizations of virtual observations $Y_{t_p:\tau}^\ddagger$. We denote the law of $\upmu_\tau$ by $\Pi_\tau$.

The stochastic optimal control problem, defined by (\ref{eq:cost}) 
and the associated Kushner--Stratonovich equation driven by
$Y_{t_p:t_p+T}^\ddagger$, falls under the broad class of stochastic mean field control problems with a common noise driver ${\rm Y}_{t_p:t_p+T}^\ddagger$; often also referred to in this specific setting as POMDP. The desired optimal control laws ${\rm u}_\tau(\mu)$, $\tau \in [t_p,t_p+T)$ follow from the Hamilton--Jacobi--Bellman equation for the value functional $v_\tau(\mu)$ and the minimum principle \cite{Nisio15,Bensoussan92}. 

Let us assume that solving the planning problem on a computer requires $\delta t>0$ units of time. Then $t_p := t_\ast - \delta t$, $T>2\delta t$, and the planning process is repeated every $\delta t$ units of time. The control applied to the physical process is given by $U_t = {\rm u}_t(\mu_t)$, $t \in [t_p+\delta t,t_p+2\delta t)$ with $\mu_t$ the law of $Z_t|Y_{0:t}^\dagger$.

%%%%%%%%%%%%%%%%%%%%%%%%%%%%%%
%
%%%%%%%%%%%%%%%%%%%%%%%%%%%%%%

\section{Data assimilation} \label{sec:DA}

In this section, we summarize the FPF reformulation of the
Kushner--Stratonovich equation 
\begin{equation} \label{eq:KSE}
    \partial_t \mu_t = -\nabla_z \cdot(\mu_t f)
    + \frac{1}{2}\nabla_z \cdot \left(\nabla_z \cdot(\mu_t\Gamma \Gamma^{\rm T})\right) + \mu_t 
    \biggl< h-\mu_t[h],\dot{Y}_t^\dagger - \mu_t[h]\biggr>
\end{equation}
$\langle a,b\rangle = a^{\rm T} b$, for the filtering distribution $\mu_t$ as originally proposed in \cite{meyn13} and developed further in \cite{reich2019data}. 
We also discuss the EnKBF as an approximation \cite{bergemann2012ensemble,TdWMR17}.

For a given control $U_t$, the FPF formulation of \cite{reich2019data} 
assimilates the data $\dot{Y}_t^\dagger$ collected from the physical process into the virtual model using mean field differential equations
\begin{subequations} \label{eq:virtualDA}
	\begin{align}
		\dot{Z}_t &= f(Z_t,U_t) + a_t(Z_t) + \Gamma \dot{V}_t 
          + K_t(Z_t)
        \left( \dot{Y}_t^\dagger - \dot{Y}_t\right),\\
		\dot{Y}_t &= h(Z_t) + \dot{W}_t
	\end{align}
\end{subequations}
with drift term
\begin{equation}
    a_t(Z_t) = \nabla_z \cdot (K_t(Z_t) K_t(Z_t)^{\rm T}) - K_t(Z_t) \nabla_z K_t(Z_t)^{\rm T}.
\end{equation}
Let $\mu_t$ denote the law of $Z_t$. The generalized Kalman gain $K(x,\mu)$ has to satisfy the partial  differential equation 
\begin{equation} \label{eq:Poisson}
    \nabla_z \cdot (\mu K^{\rm T}) = -\mu(h-\mu[h])
\end{equation}
and $K_t(z) := K(z,\mu_t)$, which follows from the requirement that $\mu_t$ obeys the Kushner--Stratonovitch equation (\ref{eq:KSE})
\cite{meyn13,reich2019data,CRS22}.

We sketch a proof that the law of $Z_t$ indeed satisfies the Kushner--Stratonovich equation (\ref{eq:KSE}). In order to verify the formal correctness of the mean field formulation (\ref{eq:virtualDA}), we only need to analyze the data-driven contributions in (\ref{eq:virtualDA}); {\it i.e.}
\begin{subequations} \label{eq:FPP_split}
    \begin{align} \nonumber
        \dot{Z}_t =&\,\,  a_t(Z_t)  + K_t(Z_t)
        \left( \dot{Y}_t^\dagger - h(X_t)- \dot{W}_t\right)\\
        =& \,\,\frac{1}{2} a_t(Z_t)
- K_t(Z_t) \dot{W}_t   
- \frac{1}{2} K_t(Z_t) \left(h-\mu_t[h] \right) \\
&\,+\, \frac{1}{2} a_t(Z_t)
+ K_t(Z_t) \left(\dot{Y}_t^\dagger  - \mu_t[h]\right)
- \frac{1}{2} K_t(Z_t) \left(h-\mu_t[h] \right).
    \end{align}
\end{subequations}
We now investigate the time evolution of the law $\mu_t$ under the two contributions (\ref{eq:FPP_split}a) and (\ref{eq:FPP_split}b) separately. First, note that
\begin{align} \nonumber
\dot{Z}_t &= \frac{1}{2}a_t(Z_t) - K_t(Z_t) \dot{W}_t   
- \frac{1}{2} K_t(Z_t) \left(h-\mu_t[h] \right)\\ \label{eq:FPFsub}
&= - K_t(Z_t) \diamond \dot{W}_t   
- \frac{1}{2} K_t(Z_t) \left(h-\mu_t[h] \right),
\end{align}
where $\diamond$ stands for Stratonovich integration with respect to $Z_t$ and 
It\^o integration with respect to $\mu_t$ \cite{CRS22}. Hence, (\ref{eq:FPFsub})
and (\ref{eq:Poisson}) lead to the Fokker--Planck equation \cite{Pavliotis2016}
\begin{equation} \label{eq:FPE1}
    \partial_t \mu_t = \frac{1}{2}\nabla_z \cdot \left(K_t \nabla_z
    \cdot (\mu_t K_t^{\rm T} )\right) + \frac{1}{2} \nabla_z \cdot \left(\mu_t K_t (h-\mu_t[h])\right) =0.
\end{equation}
Second, we find that 
\begin{align} \nonumber
\dot{Z}_t &= \frac{1}{2} a_t(Z_t)
+ K_t(Z_t) \left(\dot{Y}_t^\dagger  - \mu_t[h]\right) 
- \frac{1}{2} K_t(Z_t) \left(h-\mu_t[h] \right)\\
&= K_t(Z_t) \diamond \left(\dot{Y}_t^\dagger - \mu_t[h]\right) 
- \frac{1}{2} K_t(Z_t) \left(h-\mu_t[h] \right)
\end{align}
results in the stochastic Fokker--Planck equation
\begin{equation}\label{eq:FPE2}
    \partial_t \mu_t  
    = -\nabla_z \cdot (\mu_t K_t) \, \left(\dot{Y}_t^\dagger - \mu_t[h]\right)=
    \mu_t \bigg< h-\mu_t[h],\dot{Y}_t^\dagger-\mu_t[h]\biggr>
\end{equation}
provided that the quadratic variation of the common noise $Y_t^\dagger$ is equivalent to that of $d_y$-dimensional Brownian motion $W_t$. Combining these results, it follows that the law $\mu_t$ of the mean field process
\begin{align} \nonumber
\dot{Z}_t &= a_t(Z_t)  + K_t(Z_t)
        \left( \dot{Y}_t^\dagger - h(Z_t)-\dot{W}_t\right)\\
        &=K_t(Z_t) \diamond
        \left( \dot{Y}_t^\dagger - h(Z_t)-\dot{W}_t\right)
\end{align}
satisfies the Fokker--Planck equation
\begin{equation}
 \partial_t \mu_t =  \mu_t \biggl< h-\mu_t[h],\dot{Y}_t^\dagger - \mu_t[h]\biggr>.
\end{equation}
Hence, the law $\mu_t$ of $Z_t$, as induced by the mean-field evolution equation (\ref{eq:virtualDA}), satisfies the Kushner--Stratonovich equation (\ref{eq:KSE}).

We now provide an interacting particle approximation of
(\ref{eq:virtualDA}). An ensemble of $M$ particles $Z_t^{(i)}$
is propagated under the stochastic differential equations
\begin{subequations} \label{eq:fpf}
    \begin{align}
\dot{Z}_t^{(i)} =& \,\,f(Z_t^{(i)},U_t) + a_t(Z_t^{(i)}) + \Gamma
\dot{V}_t^{(i)} 
+ \hat{K}_t(Z_t^{(i)}) \left( \dot{Y}_t^\dagger - \dot{Y}_t^{(i)}\right),\\
\dot{Y}_t^{(i)} =& \,\,h(Z_t^{(i)}) + \dot{W}_t^{(i)},
    \end{align}
\end{subequations}
$i=1,\ldots,M$, with common noise $\dot{Y}_t^\dagger$ and Kalman 
gain matrix
\begin{equation} \label{eq:ideal particle gain}
    \hat{K}_t(z) := K(z,\hat{\mu}_t).
\end{equation}
Here, $\hat{\mu}_t$ denotes the empirical filtering distribution 
\begin{equation} \label{eq:particle empirical measure}
    \hat{\mu}_t({\rm d}z) = \frac{1}{M}\sum_{j=1}^M
    \delta(z-Z_t^{(j)}),
\end{equation}
$\delta(\cdot)$ the Dirac measure, induced by the particles $\{Z_t^{(j)}\}$ at time $t$.

\begin{example} \label{example1}
In practice, the Kalman gain matrix $K(z,\mu)$, which satisfies (\ref{eq:Poisson}), is not available and further approximations are required in (\ref{eq:ideal particle gain}) \cite{TdWMR17}. A particularly popular approximation is provided by the so called constant gain approximation
\begin{equation} \label{eq:KB gain}
K_t = C_t^{zh},
\end{equation}
where $C_t^{zh}$ denotes the covariance matrix between $Z_t$ and $h(Z_t)$ according to the law $\mu_t$. Since 
(\ref{eq:KB gain}) does not explicitly depend on the state $Z_t$, the additional drift term $a_t(z)$ vanishes.
The constant gain approximation turns the interacting particle system (\ref{eq:fpf}) into the EnKBF, as first proposed in \cite{bergemann2012ensemble}, with the particle based Kalman gain defined by
\begin{align} \nonumber
    \hat{K}_t &= \frac{1}{M-1}\sum_{j=1}^M \left(Z_t^{(j)}-\frac{1}{M}
    \sum_{i=1}^M Z_t^{(i)}\right)
    h(Z_t^{(j)})^{\rm T} \\
    &= \frac{1}{M-1}\sum_{j=1}^M Z_t^{(j)}
    \left( h(Z_t^{(j)}) -\frac{1}{M}
    \sum_{i=1}^M h(Z_t^{(i)})\right)^{\rm T}.
    \label{EnKBF gain}
\end{align}
An even greater simplification arises from a truly constant gain $K$ in which case we can set $M=1$
in (\ref{eq:fpf}) to obtain the stochastic differential equations
\begin{subequations} \label{eq:fpf nudging}
    \begin{align}
\dot{Z}_t =& \,\,f(Z_t,U_t)  + \Gamma
\dot{V}_t + K \left( \dot{Y}_t^\dagger - \dot{Y}_t\right),\\
\dot{Y}_t =& \,\,h(Z_t) + \dot{W}_t.
    \end{align}
\end{subequations}
\end{example}
 
%%%%%%%%%%%%%%%%%%%%%%%%%%%%%%
%
%%%%%%%%%%%%%%%%%%%%%%%%%%%%%%

\section{Mean field PMP for stochastic optimal control} \label{sec:SOC}

%%%%%%%%%%%%%%%%%%%%%%%%%%%%%%
%
%%%%%%%%%%%%%%%%%%%%%%%%%%%%%%

In this section, we summarize the mean field  PMP formulation of \cite{R26a} for fully observed diffusion processes and generalize it to state dependent diffusion matrices $\Sigma(x) \in \mathbb{R}^{d_x\times d_x}$. More specifically, we consider fully observed controlled diffusion processes of the form
\begin{equation} \label{eq:fSDE}
\dot{X}_t = b(X_t,U_t) + Q(X_t) \dot{B}_t
\end{equation}
for given drift $b:\mathbb{R}^{d_x}\to \mathbb{R}^{d_x}$, 
state dependent diffusion matrix $Q(x)\in \mathbb{R}^{d_x\times d_b}$ such that $\Sigma(x) = Q(x)Q(x)^{\rm T}$, 
and initial condition $X_0 = \bar x_0$. 
Here, the dimension $d_x$ of the state variable $X_t$ and the dimension $d_b$ of Brownian motion $B_t$ may be different. The cost function
for given control $U_{0:T}$ is assumed to be of the form
\begin{equation} \label{eq:cost_standard}
J(U_{0:T};\bar x_0) = \mathbb{E} \left[
\int_0^T \left(c(X_t) + \frac{1}{2}\|U_t\|_R^2  \right){\rm d}t + q(X_T) \right]
\end{equation}
with expectation taken over all realizations of $X_{0:T}$ given
the initial condition $X_0 = \bar x_0$,  The running cost $c(x)$ and the terminal cost $q(x)$ are given. 

It is assumed that the Hamiltonian 
\begin{equation}
    \mathcal{V}(x,p,u) = \bigl< b(x,u),p \bigr> + \frac{1}{2}\|u\|_R^2 + c(x)
\end{equation}
of the classical PMP \cite{pontryagin,Bensoussan} has a unique minimizer $u_\ast$ for all pairs $(x,p)$; {\it i.e.}
\begin{equation}
    \nabla_u \bigl< b(x,u_\ast),p\bigr> + R^{-1} u_\ast = 0,
\end{equation}
which gives rise to the function $u_\ast (x,p)$. Given the value function $v_t(x)$, which satisfies the associated Hamilton--Jacobi--Bellman equation \cite{Bensoussan}
\begin{equation} \label{eq:HJB}
    -\partial_t v_t(x) = \min_u \mathcal{V}(x,\nabla_x v_t(x),u)   
    + \frac{1}{2} \Sigma(x) : D_x^2 v_t(x), \qquad v_T(x) = q(x),
\end{equation}
with the Frobenius inner product $A:B = \mbox{tr}\,(A B^{\rm T})$, 
the desired feedback control law becomes
\begin{equation} \label{eq:feedback control law}
    u_t(x) = u_\ast (x,\nabla_x v_t(x)).
\end{equation}
Instead of solving the Hamilton--Jacobi--Bellman equation numerically, the mean field PMP formulation of \cite{R26a} leads to an interacting particle approximation, which we will combine in Section \ref{sec:POMDP} with the FPF formulation (\ref{eq:fpf}) from Section \ref{sec:DA}. We assume throughout this paper that the value function $v_t$, as determined by (\ref{eq:HJB}), is sufficiently smooth.

\begin{example} \label{ex:LQR}
Linear quadratic regulator (LQR) problems amount to the specfic form
\begin{equation}
    b(x,u) = b_0(x) + G(x)u
\end{equation}
in which case the function $u_\ast(x,p)$ is given by 
\begin{equation}
u_\ast(x,p) = -R G(x)^{\rm T}p
\end{equation}
and the optimal feedback control law (\ref{eq:feedback control law}) becomes
\begin{equation}
    u_t(x) = -R G(x)^{\rm T}\nabla_x v_t(x).
\end{equation}
\end{example}

We now summarize the mean field PMP formulation of \cite{R26a}. In the mean field approach,  states $X_t$ and co-states $P_t$ become functions of labels $a\in \Omega = [0,1]^{d_x}$; {\it i.e.}, $X_t(a)$ and $P_t(a)$. 
We assume a uniform distribution on label space $\Omega$ and choose $X_0:\Omega \to \mathbb{R}_{d_x}$ such that $X_0 \sim {\rm N}(\bar x_0,C_0^{xx})$ for an appropriately chosen covariance matrix $C_0^{xx}$. The terminal conditions are $P_T(a) = \nabla_x q(X_T(a))$. The distribution of $X_t(\cdot)$ at time $t\ge 0$ is denoted by $\rho_t(x)$.

We introduce the functions $\beta_t(a)\in \mathbb{R}^{d_x}$ and $\phi_t(x) \in \mathbb{R}$. In fact, $\phi_t(x)$ corresponds to the value function $v_t(x)$ and is defined implicitly through the condition
\begin{equation} \label{eq:constraint}
P_t(a) = \nabla_x \phi_t(X_t(a))
\end{equation}
for all $a\in \Omega$. The function $\beta_t(a)$ arises as the Lagrange multiplier for the constraint (\ref{eq:constraint}).

The Hamiltonian for the mean field control formulation with state dependent $\Sigma(x)$ is now given by
\begin{align} \nonumber
\mathrm{H}(P_t,X_t,U_t,\beta_t,\phi_t) =&
\int_\Omega \biggl\{ 
\mathcal{V}(X_t(a),P_t(a),U_t(a)) +  \frac{1}{2} \Sigma (X_t(a)) : D_x^2 \phi_t(X_t(a))
%\bigl< 
%b(X_t(a),U_t(a)),P_t(a) \bigr> \biggr. +c(X_t(a)) + \frac{1}{2}\|U_t(a)\|_R^2 
\\ 
&\biggl. \quad  \,+ \,\biggl< \beta_t(a), \nabla_x \phi_t(X_t(a))-P_t(a)\biggr> 
\biggr\}{\rm d}a. \label{eq:Hamiltonian_function_Lagrange}
\end{align}
The associated constrained Hamiltonian equations of motion become
\begin{subequations} \label{eq:Hamiltonian_Lagrange}
    \begin{align}
    \dot{X}_t &= \frac{\updelta \mathrm{H}}{\updelta P_t} =
    b(X_t,U_t) - \beta_t ,\\ \nonumber
        -\dot{P}_t &= \frac{\updelta \mathrm{H}}{\updelta X_t} = \left(D_x b(X_t,U_t)\right)^{\rm T}P_t + \nabla_x c(X_t) \\
        & \quad \qquad \quad \,\,\,\,+ \,\frac{1}{2}\nabla_x (\Sigma:D_x^2 \phi_t)(X_t) + D_x^2 \phi_t(X_t)\,\beta_t,\\
        0 &= \frac{\updelta \mathrm{H}}{\updelta \beta_t} =  \nabla_x \phi_t(X_t) -P_t,\\
        0 &= \frac{\updelta \mathrm{H}}{\updelta \phi_t} = -\nabla_x\cdot\left(
        \rho_t \left(\beta_t - \frac{1}{2} \nabla_x\cdot 
        \Sigma\,(X_t)- \frac{1}{2}\Sigma\nabla_x \log\rho_t\,(X_t)\right)\right),\\
        0&= \frac{\updelta \mathrm{H}}{\updelta U_t}=   
        \nabla_u \langle b(X_t,U_t),P_t\rangle + R^{-1} U_t.
    \end{align}
\end{subequations}
Therefore, condition (\ref{eq:constraint}) arises as a constraint and a natural choice for $\beta_t$ is
\begin{equation} \label{eq:beta natural}
\beta_t^\ast(a) = \frac{1}{2} \nabla_x\cdot \Sigma (X_t(a))+ \frac{1}{2}\Sigma(X_t(a))\, 
\nabla_x \log\rho_t(X_t(a)).
\end{equation}

A key observation from \cite{R26a} is that the time evolution of $\nabla_x \phi_t$ is independent of the specific choice of $\beta_t$. Indeed, differentiating the constraint (\ref{eq:Hamiltonian_Lagrange}c) with respect to time results in
\begin{equation} \label{eq:verification}
    \partial_t \nabla_x  \phi_t\,(X_t) = \dot{P}_t -
    D_x^2 \phi_t(X_t)\dot{X}_t = \dot{P}_t +
    D_x^2 \phi_t(X_t)\left(\beta_t-b(X_t,U_t)\right)
\end{equation}
and using (\ref{eq:Hamiltonian_Lagrange}b) show that the time evolution of $\nabla_x \phi_t$ is independent of the chosen $\beta_t$. Furthermore, the resulting evolution equation agrees with that of the gradient $\nabla_x v_t$ of the value function $v_t$, which satisfies the Hamilton--Jacobi--Bellman equation (\ref{eq:HJB}) \cite{R26a}.

This observation allows us to decouple the PMP boundary value problem. Given a reference 
control $u_t^{\rm ref}(x)$, such a decoupling is achieved, for example, through the use of
\begin{equation} \label{eq:beta}
  \beta_t(a) = b(X_t,U_t) - b(X_t,u_t^{\rm ref}(X_t)) +  \beta_t^\ast(a).
\end{equation}
Thus, (\ref{eq:Hamiltonian_Lagrange}a) can be solved forward in time without knowing $U_t$. 
Once $X_{0:T}$ is available, the remaining equations (\ref{eq:Hamiltonian_Lagrange}b)-(\ref{eq:Hamiltonian_Lagrange}c) and (\ref{eq:Hamiltonian_Lagrange}e) can be solved backwards in time for $(P_t,\nabla_x \phi_t,U_t)$. The reference control could be chosen as the optimal control from the previous receding horizon window shifted by $\delta t$ units of time. See \cite{R26a} for further details. Self consistency can be achieved by solving the forward differential equation (\ref{eq:Hamiltonian_Lagrange}a) in $X_t$ one more time with $u^{\rm ref}_t = u_t(x)$ and $u_t$ the optimal feedback control law (\ref{eq:feedback control law}).

Note that the mean field formulation proposed in \cite{R26a} assumes a constant diffusion matrix $\Sigma^{1/2}=Q$ in (\ref{eq:fSDE}). The extension (\ref{eq:Hamiltonian_Lagrange}) to non-square and state dependent $Q$ arises naturally from the Hamiltonian (\ref{eq:Hamiltonian_function_Lagrange}). The extension is essential for the receding horizon control formulation proposed in Section \ref{sec:POMDP}.

\begin{remark}
Forward and backward stochastic differential equations (FBSDEs) provide an alternative approach to stochastic optimal control problems \cite{Bensoussan,Carmona}.
In line with such a stochastic PMP we can formally set
$\beta_t = Q(X_t)\dot{B}_t$ in (\ref{eq:Hamiltonian_Lagrange}) and use the Stratonovich interpretation \cite{Pavliotis2016} for the multiplicative noise terms in (\ref{eq:verification}) and
(\ref{eq:Hamiltonian_Lagrange}b). The time evolution of $\nabla_x \phi_t$ is then again found to remain deterministic and agrees with that of $\nabla_x v_t$ arising from the Hamilton--Jacobi--Bellman equation (\ref{eq:HJB}). In contrast to the stochastic PMP \cite{Carmona}, our mean field approach is strictly Hamiltonian and leads to deterministic evolution equations. See \cite{R25a} for an implementation of stochastic PMP in the context of digital twins.
\end{remark}

A particle based discretization of the mean field Hamiltonian system (\ref{eq:Hamiltonian_Lagrange}) gives rise to $L$ interacting particles $(X_t^{(l)},P_t^{(l)})$
with evolution equations
\begin{subequations} \label{eq:PMP particles}
    \begin{align}
    \dot{X}_t^{(l)} &= b(X_t^{(l)},U_t^{(l)}) - \beta_t^{(l)} ,\\ \nonumber
        -\dot{P}_t^{(l)} &=  \left(D_x b(X_t^{(l)},U_t^{(l)})\right)^{\rm T}P_t^{(l)} + \nabla_x c(X_t^{(l)}) \\
        & \quad \,\,+ \,\frac{1}{2}\nabla_x (\Sigma:D_x \hat \psi_t)(X_t^{(l)}) + D_x \hat \psi_t(X_t^{(l)})\,\beta_t^{(l)},
    \end{align}
\end{subequations}
$l=1,\ldots,L$, $P_T^{(l)} = \nabla_x q(X_T^{(l)})$, control
\begin{equation}
U_t^{(l)} = u_\ast(X_t^{(l)},P_t^{(l)}),
\end{equation}
and function $\hat \psi_t \approx \nabla_x \phi_t$, which 
is obtained from the set of particles $\{(X_t^{(l)},P_t^{(l)})\}$ by linear or non-linear regression. Nadaraya--Watson kernel regression \cite{Bierens} 
\begin{equation} \label{eq:NWKR}
\hat \psi_t(x) = \sum_{l=1}^L \omega_t^{(l)}(x) \,P_t^{(l)}, \qquad
\omega_t^{(l)}(x) = 
\frac{\exp \left(-\frac{1}{2\delta}\|x - X_t^{(l)}\|^2 \right)}{\sum_{k=1}^L \exp \left(-\frac{1}{2\delta}\|x - X_t^{(k)}\|^2 \right)},
\end{equation}
has been used, for example, in the numerical experiments of \cite{R26a}. Here $\delta >0$ is a tuning parameter which depends on $L$ and the regularity of the underlying value function $v_t(x)$.
The drift terms $\{\beta_t^{(l)}\}$ can, for example, be obtained by approximating (\ref{eq:beta}) using score estimation \cite{MRO20}. See \cite{R26a} for further details and alternative approximations.

\begin{example} \label{example2}
    We recall from \cite{R26a} the following affine approximations to the score $\nabla_x \log \rho_t(x)$ and
    $\nabla_x  \phi_t(x)$. Let $C_t^{xx}$ denote the covariance matrix of $X_t$ and $C_t^{px}$ the covariance matrix
    between $P_t$ and $X_t$; {\it i.e.}
    \begin{equation}
        C_t^{xx} = \int_\Omega (X_t-\bar x_t)(X_t-\bar x_t)\,{\rm d}a, \quad
        C_t^{px} = \int_\Omega (P_t-\bar p_t)(X_t-\bar x_t)\,{\rm d}a
    \end{equation}
    with mean $\bar x_t$ and $\bar p_t$. Then (\ref{eq:beta natural}) becomes
    \begin{equation} \label{eq:app1}
        \hat{\beta}_t^\ast(a) =   -\frac{1}{2}
        \rho_t\left[\Sigma\right] (C_t^{xx})^{-1}(X_t(a)-\bar x_t)
    \end{equation}
    and $\nabla_x \phi_t$ is approximated by
    \begin{equation} \label{eq:app2}
        \hat \psi_t(x) = C^{px}_t (C^{xx}_t)^{-1}(x-\bar x_t) + \bar p_t.
    \end{equation}
    The approximations (\ref{eq:app1})-(\ref{eq:app2}) are analogous to the constant gain approximation (\ref{eq:KB gain}) in the EnKBF. Furthermore, approximations (\ref{eq:app1}) and
    (\ref{eq:app2}) arise as critical points from the parametrized Hamiltonian
    \begin{align} \nonumber
\mathrm{H}(P_t,X_t,U_t,\mathtt{B}_t,\mathtt{c}_t,\mathtt{A}_t,\mathtt{d}_t) =&
\,\int_\Omega \biggl\{
\mathcal{V}(X_t(a),P_t(a),U_t(a)) +  \frac{1}{2} \Sigma (X_t(a)) : \mathtt{A}_t
%\bigl< 
%b(X_t(a),U_t(a)),P_t(a) \bigr> \biggr. +c(X_t(a)) + \frac{1}{2}\|U_t(a)\|_R^2 
\\ %\nonumber 
&\biggl.
%b(X_t(a),U_t(a)),P_t(a) \bigr> \biggr. +c(X_t(a)) + \frac{1}{2}\|U_t(a)\|_R^2 \\ \nonumber 
\quad \,+ \,\biggl< \mathtt{B}_tX_t(a)+\mathtt{c}_t, \mathtt{A}_tX_t(a)+\mathtt{d}_t-P_t(a)\biggr> 
\biggr\}{\rm d}a; \label{eq:Hamiltonian parametrized}
\end{align}
{\it i.e.}
\begin{subequations}
    \begin{align}
    0 = \nabla_\mathtt{A}\mathrm{H} &= \int_\Omega
    \left( (\mathtt{B}_tX_t+\mathtt{c}_t)X_t^{\rm T} + \frac{1}{2}\Sigma(X_t(a)) \right) {\rm d}a,\\
    0 = \nabla_\mathtt{d} \mathrm{H} &= \int_\Omega
    (\mathtt{B}_tX_t + \mathtt{c}_t) {\rm d}a,\\
    0 = \nabla_\mathtt{B} \mathrm{H} &= \int_\Omega
    \left(\mathtt{A}_tX_t+\mathtt{d}_t-P_t\right)X_t^{\rm T}{\rm d}a,\\
    0 = \nabla_\mathtt{c} \mathrm{H} &= \int_\Omega
    (\mathtt{A}_tX_t + \mathtt{d}_t-P_t) {\rm d}a
    \end{align}
\end{subequations}
and are therefore strictly variational. An interacting particle approximation can be derived from the discretized Hamiltonian
\begin{align} \nonumber
\mathrm{H}_L(\{P_t^{(l)}\},\{X_t^{(l)}\},\{U_t^{(l)}\},B_t,c_t,A_t,d_t) =&
\,\frac{1}{L}\sum_{l=1}^L \biggl\{ 
\mathcal{V}\left(X_t^{(l)},P_t^{(l)},U_t^{(l)}\right) + \frac{1}{2} \Sigma (X_t^{(l)}) : A_t  \biggr.   \\ %\nonumber 
&\biggl.\,+ \,\biggl< B_tX_t^{(l)}+c_t, A_tX_t^{(l)}+d_t-P_t^{(l)}\biggr> 
\biggr\} \label{eq:Hamiltonian discretized}
\end{align}
and the associated constrained Hamiltonian evolution equations
\begin{subequations} \label{eq:Hamiltonian evolution discretized}
    \begin{align}
        L^{-1}\dot{X}_t^{(l)} &= \nabla_{p^{(l)}} \mathrm{H}_L,\\
        -L^{-1}\dot{P}_t^{(l)} &= \nabla_{x^{(l)}} \mathrm{H}_L,\\
        0 &= \nabla_{u^{(l)}} \mathrm{H}_L,\\
        0 &= \nabla_{A}\mathrm{H}_L = \nabla_{B}\mathrm{H}_L\\
        0 &= \nabla_{c}\mathrm{H}_L = \nabla_{d}\mathrm{H}_L,
    \end{align}
\end{subequations}
$l=1,\ldots,L$. Here we have dropped the arguments in the Hamiltonian $\mathrm{H}_L$ for notational simplicity. The optimal set of parameters $(A_t,B_t,c_t,d_t)$ provides an empirical approximation to
the mean field parameters $(\mathtt{A}_t,\mathtt{B}_t,\mathtt{c}_t,\mathtt{d}_t)$. The prefactor of
$L^{-1}$ in (\ref{eq:Hamiltonian evolution discretized}a)-(\ref{eq:Hamiltonian evolution discretized}b) arises from the fact that the variational derivatives in (\ref{eq:Hamiltonian_Lagrange}) have been replaced by gradients with respect to $x^{(l)} \in \mathbb{R}^{d_x}$ and $p^{(l)}\in \mathbb{R}^{d_x}$.
\end{example}

%%%%%%%%%%%%%%%%%%%%%%%%%%%%%%
%
%%%%%%%%%%%%%%%%%%%%%%%%%%%%%%

\section{Control of partially observed processes} \label{sec:POMDP}

%%%%%%%%%%%%%%%%%%%%%%%%%%%%%%
%
%%%%%%%%%%%%%%%%%%%%%%%%%%%%%%

We address the POMDP planning problem by combining the FPF formulation (\ref{eq:fpf}) and
the mean field PMP formulation (\ref{eq:Hamiltonian_Lagrange}). The first step is to extend
(\ref{eq:virtualDA}) to the virtual planning process for $\tau \ge t_p$. Replacing the actual observations increments $\dot{Y}_\tau^\dagger$ by virtual ones $\dot{Y}_\tau^\ddagger$, we obtain
\begin{subequations} \label{eq:vitualPlanning}
	\begin{align} 
		\dot{Z}_\tau &= f(Z_\tau,U_\tau) 
        + a_\tau(Z_\tau) + \Gamma \dot{V}_\tau         
        + K_\tau(Z_\tau)\left(\dot{\rm Y}_\tau^\ddagger - \dot{\rm Y}_\tau \right),\\
		\dot{\rm Y}_\tau &= h({\rm X}_\tau) + \dot{W}_\tau
	\end{align}
\end{subequations}
subject to initial conditions $Z_{t_p} \sim \upmu_{t_p}$ with the initial distribution provided by the current filtering distribution; {\it i.e.}~$\upmu_{t_p}=\mu_{t_p}$. 
Note that the law $\upmu_\tau$ of $Z_\tau$, $\tau>t_p$, is a random probability measure with law denoted by $\Pi_\tau$.

Invoking the ideal twin scenario \cite{bain2008fundamentals}, the statistics of virtual observations can be represented by
\begin{equation} \label{eq:virtual_obs}
	\dot{\rm Y}_\tau^\ddagger = \upmu_\tau[h] + \dot{W}_\tau^\ddagger,
\end{equation}
where $W_\tau^\ddagger \in \mathbb{R}^{d_y}$ denotes Brownian motion. 

Thus, we have reformulated the partially observed optimal control problem with cost (\ref{eq:cost}) into a fully observed optimal control problem with respect to the evolving distributions $\upmu_\tau \sim \Pi_\tau$; called belief states in reinforcement learning \cite{Meyn}, which satisfy the Kushner--Stratonovich equation
\begin{equation} \label{eq:KSE_planning}
    \partial_\tau \upmu_\tau = -\nabla_z \cdot(\upmu_\tau f)
    + \frac{1}{2}\nabla_z \cdot \left(\nabla_z \cdot(\upmu_\tau\Gamma \Gamma^{\rm T})\right) + \upmu_\tau \bigg< h-\upmu_\tau[h]),\dot{W}_\tau^\ddagger\biggr>
\end{equation}
with initial condition $\upmu_{t_p} = \mu_{t_p}$ \cite{Bensoussan92,Bensoussan}. The mean field formulation (\ref{eq:vitualPlanning})-(\ref{eq:virtual_obs}) satisfies $Z_\tau|W^\ddagger_{t_p:\tau} \sim \upmu_\tau$. 

Therefore, we may replace the stochastic partial differential equation (\ref{eq:KSE_planning})
by the mean field stochastic differential equations (\ref{eq:vitualPlanning})-(\ref{eq:virtual_obs})
In order to implement (\ref{eq:vitualPlanning})-(\ref{eq:virtual_obs}) numerically, we follow the interacting particle approximation (\ref{eq:fpf}) to obtain:
\begin{subequations} \label{eq:IPS planning}
    \begin{align}
\dot{Z}_\tau^{(i)} =& \,f(Z_\tau^{(i)},U_\tau) + a_\tau(Z_\tau^{(i)}) + \Gamma
\dot{V}_\tau^{(i)} + \hat{K}_\tau(Z_\tau^{(i)}) \left( \dot{Y}_\tau^\ddagger - \dot{Y}_\tau^{(i)}\right),\\
\dot{Y}_\tau^{(i)} =& \,h(Z_\tau^{(i)}) + \dot{W}_\tau^{(i)},\\
\dot{Y}_\tau^\ddagger =&\,\hat{\upmu}_\tau[h] + \dot{W}_\tau^\ddagger,
    \end{align}
\end{subequations}
$i=1,\ldots,M$, and initial conditions $Z_{t_p}^{(i)}$ provided by (\ref{eq:fpf}).
Here, $\hat{\upmu}_\tau$ denotes the random empirical measure 
\begin{equation}
    \hat{\upmu}_\tau({\rm d}z) = \frac{1}{M}\sum_{i=1}^M
    \delta(z-Z_\tau^{(i)})
\end{equation}
generated by random particles $\{Z_\tau^{(i)}\}$ and
\begin{equation}
    \hat{\upmu}_\tau[h] = \frac{1}{M} \sum_{i=1}^M h(Z_\tau^{(i)}).
\end{equation}
The distribution of random empirical measures $\hat{\upmu}_\tau$, $\tau > t_p$, generated by (\ref{eq:IPS planning}) over realizations of $W_{t_p:\tau}^\ddagger$, provides an approximation of $\Pi_\tau$.

\begin{example} \label{ex:EnKBF planning}
    The EnKBF approximation (\ref{eq:KB gain}) of the Kalman--Bucy gain leads to the following specific form of
    (\ref{eq:IPS planning}):
    \begin{subequations} \label{eq:IPS planning EnKBF}
    \begin{align}
\dot{Z}_\tau^{(i)} =& \,f(Z_\tau^{(i)},U_\tau) + \Gamma
\dot{V}_\tau^{(i)} + \frac{1}{M-1} \sum_{j=1}^M Z_\tau^{(j)} \left(
h(Z_\tau^{(j)})-\hat{\upmu}_\tau[h] \right)^{\rm T}\left( \dot{Y}_\tau^\ddagger - \dot{Y}_\tau^{(i)}\right),\\
\dot{Y}_\tau^{(i)} =& \,h(Z_\tau^{(i)}) + \dot{W}_\tau^{(i)},\\
\dot{Y}_\tau^\ddagger =&\,\hat{\upmu}_\tau[h] + \dot{W}_\tau^\ddagger.
    \end{align}
\end{subequations}
We recall that the EnKBF does not assume that the underlying distributions are Gaussian \cite{CRS22}.
Instead, the EnKBF emerges as the best constant gain approximation to the Kalman gain $K(z,\mu)$ \cite{TdWMR17}.
The Kalman--Bucy gain (\ref{eq:KB gain}) still depends on the law of $Z_\tau$. A truly constant gain $K$ leads from (\ref{eq:IPS planning}) with $M=1$ to 
\begin{subequations} 
    \begin{align}
\dot{Z}_\tau =& \,\,f(Z_\tau,U_\tau)  + \Gamma
\dot{V}_\tau + K \left( \dot{Y}_\tau^\ddagger - \dot{Y}_\tau\right),\\
\dot{Y}_\tau =& \,\,h(Z_\tau) + \dot{W}_\tau,\\
\dot{Y}_\tau^\ddagger =&\,\,h(Z_\tau) + \dot{W}_\tau^\ddagger
    \end{align}
\end{subequations}
or, in compact form,
\begin{equation}
    \dot{Z}_\tau = f(Z_\tau,U_\tau) + \Gamma \dot{V}_\tau + \sqrt{2} K \dot{W}_\tau
\end{equation}
to which mean field PMP can be applied directly. In all other cases ($M>1$), an extended variable approach has to be followed, as explained next.
\end{example}

We introduce the extended state variable 
\begin{equation}
x = (z_1^{\rm T},z_2^{\rm T},\ldots,z_M^{\rm T})^{\rm T} \in \mathbb{R}^{d_z M};
\end{equation}
{\it i.e.}, $d_x = d_z\,M$ and
\begin{equation}
X_\tau = \left( \begin{array}{c} Z_\tau^{(1)}\\ \vdots \\ Z_\tau^{(M)}
\end{array} \right) \in \mathbb{R}^{d_zM},
\end{equation}
and rewrite (\ref{eq:IPS planning}) as a controlled diffusion process as discussed in Section \ref{sec:SOC}; {\it i.e.},
\begin{equation} \label{eq:POMDP_SDE}
\dot{X}_\tau = b(X_\tau,U_\tau) + Q(X_\tau) \dot{B}_\tau
\end{equation}
with initial condition $X_{t_p} = \bar x_{t_p}$ determined by the empirical filtering distribution 
\begin{equation}
\hat{\mu}_{t_p}({\rm d}z) = \frac{1}{M} \sum_{i=1}^M \delta(z-Z_{t_p}^{(i)}),
\end{equation}
which is provided by the FPF (\ref{eq:fpf}).
Here $B_\tau$ is the collection of Brownian motions $\{V_\tau^{(j)}\}$, $\{W_\tau^{(j)}\}$, and $W_\tau^\ddagger$. The matrix $Q(x) \in \mathbb{R}^{d_z M\times(d_zM+d_yM + d_y)}$ is made up of entries that involve both $\Gamma$ and $\hat{K}(z_j)$; {\it e.g.}, setting $M=3$
\begin{equation} \label{eq:G}
Q (x) = \left( \begin{array}{ccccccc}   \Gamma & 0 & 0 & \hat{K}(z_1) & 0 & 0 & 
\hat{K}(z_1) \\
0 & \Gamma  & 0 & 0 & \hat{K}(z_2) & 0  & \hat{K}(z_2) \\
0 & 0 & \Gamma & 0 & 0 & \hat{K}(z_3) & \hat{K}(z_3) 
\end{array} \right).
\end{equation}
The drift function
$b(x,u):\mathbb{R}^{d_z M}\to \mathbb{R}^{d_z M}$ collects the drift terms in (\ref{eq:IPS planning}); {\it i.e.}
\begin{equation} \label{eq:b}
b(x,u) = \left( \begin{array}{c} 
f(z_1,u) 
+ a_t(z_1) 
- \hat{K}(z_1)\left(h(z_1)-\hat{\upmu}[h]\right)\\
f(z_2,u) 
+ a_t(z_2) 
- \hat{K}(z_2)\left(h(z_2)-\hat{\upmu}[h]\right)\\
\vdots \\
f(z_M,u) 
+ a_t(z_M)
- \hat{K}(z_M)\left(h(z_M)-\hat{\upmu}[h]\right)
\end{array}\right).
\end{equation}
Recall that $\hat{K}(z)$ depends on all entries in $x$ through the empirical measure $\hat{\upmu}$ induced by $\{z_i\}$. 
See Section \ref{sec:Example} for more details in the context of the EnKBF.

Substituting $\upmu_\tau$ with its empirical approximation $\hat{\upmu}_\tau$, the cost function (\ref{eq:cost}) becomes
\begin{align} \nonumber
	J(U_{t_p:t_p+T};X_{t_p}) &= \mathbb{E} \left[ \int_{t_p}^{t_p+T}
	\left( \hat{\upmu}_\tau[r] + \frac{1}{2} \|U_\tau\|_R^2\right) {\rm d}\tau + \hat{\upmu}_{t_p+T}[\psi] \right]\\
    &= \mathbb{E} \left[ \int_{t_p}^{t_p+T}
	\left( \frac{1}{M}\sum_{i=1}^M r(Z_\tau^{(i)}) + \frac{1}{2} \|U_\tau\|_R^2\right) {\rm d}\tau + \frac{1}{M} \sum_{i=1}^M 
    g(Z^{(i)}_{t_p+T}) \right] \label{eq:POMDP_cost}
\end{align}
with expectation taken over the Brownian motion $B_{t_p:t_p+T}$. Hence, 
the running cost in (\ref{eq:cost_standard}) becomes
\begin{equation}
c(x) = \frac{1}{M}\sum_{i=1}^M r(z_i)
\end{equation}
and the terminal cost
\begin{equation}
q(x) = \frac{1}{M} \sum_{i=1}^M g(z_i) .
\end{equation}

Thus, we have reformulated the planning problem for partially observed diffusion processes into a standard stochastic optimal control problem, and the mean field PMP formulation (\ref{eq:Hamiltonian_Lagrange}) from Section \ref{sec:SOC} can be applied with mean field Hamiltonian
\begin{align} \nonumber
\mathcal{H}(P_\tau,X_\tau,U_\tau,\beta_\tau,\phi_\tau) =&
\int_\Omega \biggl\{ \mathcal{V}_M\left(X_\tau(a),P_\tau(a),U_\tau(a)\right)+  \frac{1}{2} \Sigma (X_\tau(a)) : D_x^2 \phi_\tau(X_\tau(a))\biggl. \\ 
&\biggl. \quad \,+ \,\biggl< \beta_\tau(a), \nabla_x \phi_\tau(X_\tau(a))-P_\tau(a)\biggr> 
\biggr\}{\rm d}a, \label{eq:POMDP Hamiltonian}
\end{align}
classical PMP Hamiltonian
\begin{equation} \label{eq:classical PMP rescalled}
    \mathcal{V}_M(x,p,u) = \bigl< 
b(x,u),p \bigr> \biggr. +M \left\{c(x) + \frac{1}{2}\|u\|_{R}^2\right\},
\end{equation}
and diffusion matrix $\Sigma(x) = Q(x)Q(x)^{\rm T}$. Note that the running and control cost has been multiplied by the size of the ensemble, $M$, in the classical PMP Hamiltonian (\ref{eq:classical PMP rescalled}). Such a rescaling does not affect the resulting optimal control; but it leads to more natural constrained Hamiltonian evolution equations; {\it e.g.}
\begin{align}\nonumber 
    0 &= 
    \nabla_u \big< b(X_\tau,U_\tau),\nabla_{x}\phi_\tau(X_\tau)\big> + M R^{-1} U_\tau 
    \\ \label{eq:implicit control law}
    &=
    \frac{1}{M}\sum_{i=1}^M \nabla_u \big< f(Z_\tau^{(i)},U_\tau),\nabla_{z_i}\phi_\tau(X_\tau)\big> + R^{-1} U_\tau,
\end{align}
and terminal condition $P_{t_p+T} = \nabla_x g(X_{t_p+T})$. The initial conditions are assumed to satisfy $X_{t_p} \sim {\rm N}(\bar x_{t_p},C_{t_p}^{xx})$ for an appropriate covariance matrix $C_{t_p}^{xx}$.

The control $U_\tau$ given by (\ref{eq:implicit control law})
is invariant under permutations of the entries $z_i$ in the vector $x = (z_1^{\rm T},\ldots,z_M^{\rm T})^{\rm T}$; therefore, the control depends only on the associated empirical measure 
\begin{equation}
\hat{\mu}({\rm d}z) = \frac{1}{M} \sum_{i=1}^M \delta (z-z_i)
\end{equation}
and provides an approximation of the optimal control law ${\rm u}_\tau(\mu)$. For the same reason, the value function $\phi_\tau(x)$ approximates the value functional $\upphi_\tau(\mu)$ \cite{Nisio15}. We denote the L-derivative \cite{CarmonaDelarue} of $\upphi_\tau(\mu)$  by $\uppsi_\tau(z,\mu) \in \mathbb{R}^{d_z}$. Then $\nabla_{z_i} \phi_\tau(x)$ provides an approximation to $\uppsi_\tau(z_i,\hat{\mu})$. The optimal control ${\rm u}_\tau(\mu)$ and the L-derivative $\uppsi_\tau(z,\mu)$ are related through 
\begin{equation}
    \mu_\tau \bigl[ \nabla_u \langle f(\,\cdot\,,{\rm u}_\tau(\mu)),\uppsi_\tau(\,\cdot\,,\mu) \rangle\bigr] + R^{-1}{\rm u}_\tau(\mu) = 0;
\end{equation}
compare (\ref{eq:implicit control law}).

Particle approximations of (\ref{eq:POMDP Hamiltonian}) should maintain the permutation invariance. For example, the standard Euclidean distance between two vectors $x = (z_1^{\rm T},\ldots,z_M^{\rm T})^{\rm T}$ and $\hat x = (\hat z_1^{\rm T},\ldots,\hat z_M^{\rm T})^{\rm T}$ should be replaced in the Nadaraya--Watson regression estimator (\ref{eq:NWKR}) by the Wasserstein distance $W_2(\hat{\mu},\tilde{\mu})$ \cite{villani2021topics} with respect to two empirical measures $\hat{\mu}$ and $\tilde{\mu}$ induced on $\mathbb{R}^{d_z}$. More precisely, let us denote the block entries in $X_\tau^{(l)}$ and $P_\tau^{(l)}$ by $\{Z_\tau^{(l,i)}\}$ and $\{\Lambda_\tau^{(l,i)}\}$, respectively.
Then, we replace (\ref{eq:NWKR}) with the regression estimator
\begin{equation} \label{eq:NWKR_2}
    \hat{\uppsi}_\tau(z,\mu) = \sum_{i=1}^M \sum_{l=1}^L
    \omega_\tau^{(i,l)}(z,\mu) \,\Lambda_\tau^{(i,l)}, 
\end{equation}
of the L-derivative $\uppsi_\tau(z,\mu)\in \mathbb{R}^{d_z}$ with normalized weights
\begin{equation}
    \omega_\tau^{(i,l)}(z,\mu) =
    \frac{\exp\left(-\frac{1}{2\delta}W_2\left(\mu,\hat{\upmu}_\tau^{(l)}\right)^2\right)
     \, \exp\left(-\frac{1}{2\delta}\|z-Z_\tau^{(l,i)}\|^2 \right)}
    {\sum_{k=1}^L \sum_{j=1}^M \exp\left(-\frac{1}{2\delta}W_2\left(\mu,\hat{\upmu}_\tau^{(k)}\right)^2\right)\, \exp\left(-\frac{1}{2\delta}\|z-Z_\tau^{(k,j)}\|^2 \right)},
\end{equation}
$\delta >0$. Alternatively, one can consider unmasked transformer architectures and build on their universal approximation properties for in-context mappings of the form $\uppsi_\tau(z,\mu) \in \mathbb{R}^{d_z}$ \cite{Transformer}. 

We recall from \cite{R26a} the freedom in choosing $\beta_\tau$, which also applies to the application of mean field PMP to (\ref{eq:POMDP_SDE}) and (\ref{eq:POMDP_cost}) and which allows for a decoupling of the associated boundary value problem. 

Once the planning process is completed and an approximation $\hat \uppsi_\tau(z,\mu)$ is available for $t_p\le \tau < t_p+T$, following (\ref{eq:implicit control law}), 
the control $U_t$ for the physical process is provided by 
\begin{equation}
0=  \hat{\mu}_t \left[\nabla_u \big\langle f(\,\cdot\,,U_t),
\hat{\uppsi}_t(\,\cdot\,,\hat{\mu}_t)\big\rangle \right]+ R^{-1} U_t
\end{equation}
with the empirical filtering distribution (\ref{eq:particle empirical measure}) obtained from the FPF formulation (\ref{eq:fpf}) using newly incoming data $\dot{Y}_t^\dagger$, $t \ge  t_p+\delta t$. Finally, $t_p \gets t_p+\delta t$ and a new planning process begins.

\begin{remark}
We contrast the proposed approach with a formulation that arises from a formal application of the separation principle \cite{Handel07}. First, we set $M = 1$; {\it i.e.}~$d_x = d_z$, and replace $\hat {\upmu}_{t_p} = \hat{\mu}_{t_p}$ by its empirical mean
\begin{equation}
    \bar x_{t_p} = \bar z_{t_p} := \hat{\mu}_{t_p}[z] = \frac{1}{M}\sum_{i=1}^M Z_{t_p}^{(i)}.
\end{equation}
Second, we reduce the planning process to (\ref{eq:virtual}a) ignoring 
(\ref{eq:virtual}b) and the impact of future observations $\dot{Y}_t^\dagger$, $t\ge t_\ast$. Hence, the planning problem is fully characterized by (\ref{eq:POMDP_SDE}) with initial condition $X_{t_p} \sim {\rm N}(\bar x_{t_p},C_{t_p}^{xx})$, $b(x,u) = f(x,u)$ and $Q(x) = \Gamma$ together with the cost function (\ref{eq:POMDP_cost}) with $c(x) = r(x)$ and $q(x) = g(x)$. Third, once $\nabla_x \phi_t(x)$ is available, the control $U_t$ for the physical process over the time interval $t \in [t_p+\delta t,t_p+2\delta t)$ satisfies
\begin{equation}
    0= \nabla_u \big\langle f(\bar z_t,U_t),\nabla_z \phi_t(\bar z_t)\big\rangle + R^{-1}U_t.
\end{equation}
Here, 
\begin{equation}
    \bar z_t = \frac{1}{M}\sum_{i=1}^M Z_t^{(i)}
\end{equation}
with $\{Z_t^{(i)}\}$ obtained from the FPF (\ref{eq:fpf}). Although this formulation is computationally less expensive than (\ref{eq:implicit control law}), it ignores uncertainties about the state ${\rm Z}_{t_p}^\dagger$ of the physical process and the impact of future observations $\dot{Y}_t^\dagger$, $t\ge t_\ast$, on the optimal control unless the physical process is linear and Gaussian; in which case the separation principle applies \cite{Bensoussan,Handel07}.
\end{remark}

%%%%%%%%%%%%%%%%%%%%%%%%%%%%
%
%%%%%%%%%%%%%%%%%%%%%%%%%%

\section{Affine approximations and EnKBF} \label{sec:Example}

In this section, we discuss the constant gain and affine approximations introduced in Examples \ref{example1}, \ref{ex:LQR}, \ref{example2}, and \ref{ex:EnKBF planning} in the context of the receding horizon planning step defined by the mean field Hamiltonian (\ref{eq:POMDP Hamiltonian}) in more detail. For given $x$, we define the Kalman--Bucy gain
\begin{equation}
    \hat{K} = \frac{1}{M-1} \sum_{i=1}^M
    \left( z_i-\frac{1}{M}\sum_{i=1}^M z_i\right)
    h(z_i)^{\rm T},
\end{equation}
which appears in $Q(x)$; compare (\ref{eq:G}), and the drift term (\ref{eq:b}). The additional drift term $a_t(z)$ vanishes under the Kalman--Bucy approximation in the limit $M\to \infty$ and is ignored. 

Affine approximations (\ref{eq:app1}) and (\ref{eq:app2}) for $\beta_t^\ast$ and
$\nabla_x \phi_t$, respectively, follow along the lines of the mean field Hamiltonian (\ref{eq:Hamiltonian parametrized}) and its interacting particle approximation (\ref{eq:Hamiltonian discretized}). 

However, propagation of chaos for the EnKBF \cite{delMoral17} implies that the entries $Z_\tau^{(i)}$ of $X_\tau$ become i.i.d.~in the mean field limit $M\to \infty$. We can exploit this property by restricting 
$A_\tau$ and $B_\tau$ to block diagonal with common ${\rm A}_\tau$ and ${\rm B}_\tau$, respectively; {\it e.g.},
\begin{equation}
    A_\tau = \left(\begin{array}{ccc} {\rm A}_\tau & 0 & 0\\
    0 & {\rm A}_\tau &0 \\ 0 & 0 & {\rm A}_\tau \end{array}\right), \qquad
    B_\tau = \left(\begin{array}{ccc} {\rm B}_\tau & 0 & 0\\
    0 & {\rm B}_\tau &0 \\ 0 & 0 & {\rm B}_\tau \end{array}\right)
\end{equation}
for $M=3$. Similarly for the vectors $c_\tau$ and $d_\tau$; {\it e.g.},
\begin{equation}
    c_\tau = \left(\begin{array}{c} {\rm c}_\tau \\
    {\rm c}_\tau  \\  {\rm c}_\tau \end{array}\right), \qquad
    d_\tau = \left(\begin{array}{c} {\rm d}_\tau \\
    {\rm d}_\tau  \\ {\rm d}_\tau \end{array}\right).
\end{equation}
Taking note of (\ref{eq:Hamiltonian discretized}), we can substitute this {\it ansatz} into 
(\ref{eq:POMDP Hamiltonian}) to obtain
\begin{align} \nonumber
&\mathcal{H}_L\left(\{P_\tau^{(l)}\},\{X_\tau^{(l)}\},\{U_\tau^{(l)}\},{\rm B}_\tau,{\rm c}_\tau,{\rm A}_\tau,{\rm d}_\tau\right) =\\
\nonumber & \qquad\qquad \qquad \qquad \qquad \qquad \qquad 
\,\frac{1}{L}\sum_{l=1}^L \biggl\{ \mathcal{V}_M\left(X_\tau^{(l)},P_\tau^{(l)},U_\tau^{(l)})\right)
%\bigl< 
%b(X_\tau^{(l)},U_\tau^{(l)}),P_\tau^{(l)} \bigr> \biggr. + M c(X_\tau^{(l)}) 
+ \frac{1}{2} \Sigma (X_\tau^{(l)}) : A_\tau \biggr. \\ 
&\qquad\qquad \qquad \qquad \qquad \qquad \qquad  \biggl.\quad \,+ \,\biggl< B_\tau X_\tau^{(l)}+c_\tau, A_\tau X_\tau^{(l)}+d_\tau-P_\tau^{(l)}\biggr> 
\biggr\} \label{eq:Hamiltonian discretized 2}
\end{align}
and the constrained Hamiltonian evolution equations
\begin{subequations}
    \begin{align}
        L^{-1}\dot{X}_\tau^{(l)} &= \nabla_{p^{(l)}} \mathcal{H}_L,\\
        -L^{-1}\dot{P}_\tau^{(l)} &= \nabla_{x^{(l)}} \mathcal{H}_L,\\
        0 &= \nabla_{u^{(l)}} \mathcal{H}_L,\\
        0 &= \nabla_{{\rm A}}\mathcal{H}_L = \nabla_{{\rm B}}\mathcal{H}_L\\
        0 &= \nabla_{{\rm c}}\mathcal{H}_L = \nabla_{{\rm d}}\mathcal{H}_L,
    \end{align}
\end{subequations}
$l=1,\ldots,L$. Here we have again dropped the arguments in the Hamiltonian $\mathcal{H}_L$ for simplicity of notation. Then
\begin{equation}
    {\rm c}_\tau = -\frac{1}{LM}\sum_{l=1}^L\sum_{i=1}^M
    {\rm B}_\tau Z_\tau^{(l,i)}, \qquad {\rm d}_\tau = \frac{1}{LM}\sum_{l=1}^L \sum_{i=1}^M
    \left(\Lambda_\tau^{(l,i)}-{\rm A}_\tau Z_\tau^{(l,i)}\right)
\end{equation}
and
\begin{equation}
{\rm A}_\tau = C_\tau^{\lambda z}(C_\tau^{zz})^{-1}, \qquad
{\rm B}_\tau = -\frac{1}{2} S_\tau (C_\tau^{zz})^{-1}
\end{equation}
with 
\begin{subequations}
    \begin{align}
    C_\tau^{zz} &= \frac{1}{LM}\sum_{i=1}^M\sum_{l=1}^L
    Z_\tau^{(l,i)} \left( Z_\tau^{(l,i)}-\frac{1}{LM}\sum_{j=1}^M\sum_{k=1}^L
    Z_\tau^{(k,j)} \right)^{\rm T}, \\
    C_\tau^{\lambda z} &= \frac{1}{LM}\sum_{i=1}^M\sum_{l=1}^L
    \Lambda_\tau^{(l,i)} \left( Z_\tau^{(l,i)}-\frac{1}{LM}\sum_{j=1}^M\sum_{k=1}^L
    Z_\tau^{(k,j)} \right)^{\rm T},
    \end{align}
\end{subequations}
and $S_\tau \in \mathbb{R}^{d_z\times d_z}$ defined as the average over the diagonal $d_z\times d_z$ blocks of $\Sigma(X_\tau^{(l)}) \in \mathbb{R}^{d_x\times d_x}$ and then again averaged over $l=1,\ldots,L$.

Let us return to the LQR problem from Example \ref{ex:LQR}. We assume that the drift function $f(x,u)$ is of the form
\begin{equation}
f(x,u) = f_0(x) + F(x)u.
\end{equation}
Then 
\begin{equation}
b_0(x) = \left( \begin{array}{c}
f_0(z_1)-\hat{K}(h(z_1)-\hat{\upmu}[h])\\
f_0(z_2) -\hat{K}(h(z_2)-\hat{\upmu}[h])\\ \vdots \\
f_0(z_M) -\hat{K}(h(z_M)-\hat{\upmu}[h]) \end{array}
\right) \in \mathbb{R}^{d_z M}
\end{equation}
and
\begin{equation}
    G(x) = \left( %\begin{array}{c}
    F(z_1)^{\rm T},F(z_2)^{\rm T}, \ldots ,F(z_M)^{\rm T}
    %\end{array}
    \right)^{\rm T} \in \mathbb{R}^{d_z M\times d_u}
\end{equation}
and the mean field optimal control becomes
\begin{equation}
u_t(x) = \frac{-1}{M}\sum_{i=1}^M 
R F(z_i)^{\rm T} \nabla_{z_i}\phi_t(x) \approx 
\frac{-1}{M} \sum_{i=1}^M R F(z_i)^{\rm T} \left( {\rm A}_t z_i + {\rm d}_t \right).
\end{equation}
The control applied to the physical process is finally given by 
\begin{equation}
    U_t = \frac{-1}{M} \sum_{i=1}^M R F(Z_t^{(i)})^{\rm T} \left( {\rm A}_t Z_t^{(i)} + {\rm d}_t \right)
\end{equation}
with $\{Z_t^{(i)}\}$ provided by the EnKBF implementation of (\ref{eq:fpf})
for $t \in [t_p+\delta t,t_p + 2\delta t)$.

%%%%%%%%%%%%%%%%%%%%%%%%%%%%%%
%
%%%%%%%%%%%%%%%%%%%%%%%%%%%%%%

\section{Concluding remarks} \label{sec:conclusions}

We have outlined a mathematical framework for digital twins that are updated continuously in time. The proposed interacting particle formulation bypasses the separation of data assimilation and control tasks, which can be justified in a linear Gaussian setting; but is known to be suboptimal otherwise \cite{Handel07,PGM}.

Both, the assimilation of data and the virtual planning process rely on the ideal twin assumption which is not satisfied in practical applications. This raises the question of how to assess the quality, robustness, and reliability of a particular computational implementation. A natural approach to compare the performance of different computational models of one and the same physical process is to score them  with respect to a data dependent quantity of interest (QoI) $w(y)$ \cite{gneiting07}. More precisely, the continuously ranked probability score allows for a ranking based on predicted QoI $\{w_t^{(i)} = w(Y_t^{(i)})\}$ and actual observed $w_t^\dagger = w(Y_t^\dagger)$ \cite{broecker12,reich2015probabilistic}. A rigorous statistical investigation under model misspecification is left for future research.

We have not discussed numerical implementations of the proposed mathematical framework beyond a basic linear approximation {\it ansatz}. Scalability with respect to dimension $d_x = d_z M$ of the planning process $X_t$ will be crucial. Application of the EnKBF should allow for moderate ensemble sizes $M$ even for relatively high dimensional state spaces $Z_t \in \mathbb{R}^{d_z}$. Another crucial issue is the effective approximation of $\nabla_x \phi_t(x)$ under permutation invariance and suitable choices for $\beta_t(a)$ beyond the affine approximation introduced in Section \ref{sec:Example}. As briefly discussed in Section \ref{sec:POMDP}, unmasked transformers \cite{Transformer} and Nadaraya--Watson regression of the form (\ref{eq:NWKR_2}) provide possible starting points. 

The proposed formulation should be compared numerically with, for example, output feedback MPC \cite{MPC,Allgower}, belief state covariance steering \cite{9811560}, and forward-backward SDE optimal control formulations \cite{Carmona,R25a}.

Many practical applications of the digital twin paradigm will lead to discrete in time observations $Y_{t_n}^\dagger$ with $t_n = n\,\Delta t$, $n\ge 1$, and $\Delta t >0$. An extension of the proposed framework could follow the mollification approach of \cite{br10b}, which reformulates sequential data assimilation as a continuous in time mean field ODE system to which the control framework of Section \ref{sec:SOC} is applicable. 

\medskip

\noindent
{\bf Acknowledgment.} This work has been funded by Deutsche Forschungsgemeinschaft (DFG) - Project-ID 318763901 - SFB1294. Part of the work was initiated while the author was visiting the Institute for Mathematical Sciences, National University of Singapore, in 2026.

\bibliographystyle{plainurl}
%%%%%%%%%%%%%%%%%%%%%%%%%%%%%%%%%
%
\bibliography{bib-database}

\begin{thebibliography}{10}

\bibitem{bain2008fundamentals}
Alan Bain and Dan Crisan.
\newblock {\em Fundamentals of Stochastic Filtering}, volume~60.
\newblock Springer Science \& Business Media, 2008.

\bibitem{BNST25}
Christian Bayer, Saifeddine~Ben Naamia, Erik von Schwerin, and Raul Tempone.
\newblock A {P}ontryagin maximum principle on the belief space for
  continuous-time optimal control with discrete observations.
\newblock Technical report, arXiv:2512.24916, 2025.

\bibitem{Bensoussan92}
Alain Bensoussan.
\newblock {\em Stochastic Control of Partially Observable Systems}.
\newblock Cambridge University Press, Cambridge, 1992.

\bibitem{Bensoussan}
Alain Bensoussan.
\newblock {\em Estimation and Control of Dynamical Systems}.
\newblock Springer, Cham, 2018.

\bibitem{br10b}
Kay Bergemann and Sebastian Reich.
\newblock A mollified ensemble {K}alman filter.
\newblock {\em Q. J. R. Meteorological Soc.}, 136:1636--1643, 2010.

\bibitem{bergemann2012ensemble}
Kay Bergemann and Sebastian Reich.
\newblock An ensemble {K}alman-{B}ucy filter for continuous data assimilation.
\newblock {\em Meteorologische Zeitschrift}, 21(3):213, 2012.

\bibitem{Bierens}
Herman~J. Bierens.
\newblock The {N}adaraya–-{W}atson kernel regression function estimator.
\newblock In {\em Topics in Advanced Econometrics}, page 212–247. Cambridge
  University Press, New York, 1994.

\bibitem{broecker12}
Jochen Br\"ocker.
\newblock Evaluating raw ensembles with the continuous ranked probability
  score.
\newblock {\em Quarterly J. Royal Meteor. Soc.}, 138:1611--1617, 2012.

\bibitem{CRS22}
Edoardo Calvello, Sebastian Reich, and Andrew~M. Stuart.
\newblock Ensemble {K}alman methods: {A} mean field perspective.
\newblock {\em Acta Numerica}, 34:123--291, 2025.

\bibitem{Carmona}
Ren\'e Carmona.
\newblock {\em Lectures on {BSDE}s, Stochastic Control, and Stochastic
  Differential Games with Financial Applications}.
\newblock SIAM, Philadelphia, 2016.

\bibitem{CarmonaDelarue}
Ren\'e Carmona and Francois Delarue.
\newblock {\em Probabilistic Theory of Mean Field Games with Applications {I}}.
\newblock Springer-Verlag, Cham, 2018.

\bibitem{Evensenetal2022}
Geir Evensen, Femke~C. Vossepoel, and Peter~Jan. {van Leeuwen}.
\newblock {\em Data Assimilation Fundamentals: {A} unified Formulation of the
  State and Parameter Estimation Problem}.
\newblock Springer Nature Switzerland AG, Cham, Switzerland, 2022.

\bibitem{Allgower}
Rolf Findeisen and Frank Allg\"ower.
\newblock An introduction to nonlinear model predictive control.
\newblock In A.~G. Jager, de and H.~J. Zwart, editors, {\em Systems and control
  : 21th Benelux meeting 2002}, pages 119--141. Technische Universiteit
  Eindhoven, 2002.

\bibitem{Transformer}
Takashi Furuya, Maarten~V. de~Hoop, and Gabriel Peyr\'e.
\newblock Transformers are universal in-context learners.
\newblock Technical report, arXiv:2408.01367.24916, 2024.

\bibitem{gneiting07}
Tilmann Gneiting, Fadoua Balabdaoui, and Adrian~E. Raftery.
\newblock Probabilistic forecasts, calibration and sharpness.
\newblock {\em J.~R.~Statist.~Soc.~B}, 69:243--268, 2007.

\bibitem{Willcox24}
Luwen Huang, Michael Kapteyn, and Karen~E. Willcox.
\newblock Digital twin: {G}raph formulations for managing complexity and
  uncertainty.
\newblock In {\em 2024 IEEE Smart World Congress (SWC)}, pages 2141--2148,
  2024.

\bibitem{jazwinski2007stochastic}
Andrew~H Jazwinski.
\newblock {\em Stochastic Processes and Filtering Theory}.
\newblock Courier Corporation, 2007.

\bibitem{PGM}
Michael~G. Kapteyn, Jacob V.~R. Pretorius, and Karen~E. Willcox.
\newblock A probabilistic graphical model foundation for enabling predictive
  digital twins at scale.
\newblock {\em Nature Computational Science}, 1:337--347, 2021.

\bibitem{Laplace}
Pierre-Simon Laplace.
\newblock {\em A philosophical essay on probabilities}.
\newblock John Wiley and Sons, New York, 1902.

\bibitem{MRO20}
Dimitra Maoutsa, Sebastian Reich, and Manfred Opper.
\newblock Interacting particle solutions of {F}okker--{P}lanck equations
  through gradient-log-density estimation.
\newblock {\em Entropy}, 22(8), 2020.

\bibitem{PM13}
Prashant~G. Mehta and Sean~P. Meyn.
\newblock A feedback particle filter-based approach to optimal control with
  partial observations.
\newblock In {\em 52nd IEEE Conference on Decision and Control}, pages
  3121--3127, 2013.

\bibitem{Meyn}
Sean Meyn.
\newblock {\em Control Systems and Reinforcement Learning}.
\newblock Cambridge University Press, Cambridge, 2022.

\bibitem{delMoral17}
Pierre~Del Moral, Aline Kurtzmann, and Julian Tugaut.
\newblock On the stability and the uniform propagation of chaos of a class of
  extended ensemble {K}alman--{B}ucy filters.
\newblock {\em SIAM Journal on Control and Optimization}, 55:119--155, 2017.

\bibitem{NA24}
{National Academy of Sciences}.
\newblock Foundational research gaps and future directions for digital twins.
\newblock Technical report, Washington, DC: The National Academies Press, 2024.

\bibitem{Nisio15}
Makiko Nisio.
\newblock {\em Stochastic Control Theory}.
\newblock Springer, Tokyo, 2nd edition, 2015.

\bibitem{R25c}
Manfred Opper and Sebastian Reich.
\newblock Digital twins: {M}c{K}ean--{P}ontryagin control for partially
  observed physical twins.
\newblock {\em Journal Computational Physics}, 563, 2026.

\bibitem{R26a}
Manfred Opper and Sebastian Reich.
\newblock On a mean-field {P}ontryagin minimum principle for stochastic optimal
  control.
\newblock {\em SISC}, in press, 2026.

\bibitem{Pavliotis2016}
Grigorios~A. Pavliotis.
\newblock {\em Stochastic Processes and Applications}.
\newblock Springer Verlag, New York, 2016.

\bibitem{pontryagin}
Lev~S. Pointryagin, Vladimir~G. Boltyanskii, Revaz~V. Gamkrelidze, and
  Evgenii~F. Mihchenko.
\newblock {\em The Mathematical Theory of Optimal Processes}.
\newblock John Wiley \& Sons, New York, 1962.

\bibitem{ASTROM1965}
Karl~Johan \r{A}str\"om.
\newblock Optimal control of {M}arkov processes with incomplete state
  information.
\newblock {\em Journal of Mathematical Analysis and Applications}, 10:174--205,
  1965.

\bibitem{MPC}
James~B. Rawlings, David~Q. Mayne, and Moritz~M. Diehl.
\newblock {\em Model Predictive Control: {T}heory, Computation, and Design}.
\newblock Nob Hill Publishing, Madison, 2nd edition, 2018.

\bibitem{reich2019data}
Sebastian Reich.
\newblock Data assimilation: {T}he {S}chr{\"o}dinger perspective.
\newblock {\em Acta Numerica}, 28:635--711, 2019.

\bibitem{R25a}
Sebastian Reich.
\newblock {Ensemble Kalman-Bucy filtering for nonlinear model predictive
  control}.
\newblock In Barry Koren and Dirk Roose, editors, {\em Numerical Analysis and
  Scientific Computing: An Anthology}, volume 156 of {\em Lecture Notes in
  Computational Science and Engineering}, pages 309--325. Springer, 2026.

\bibitem{reich2015probabilistic}
Sebastian Reich and Colin Cotter.
\newblock {\em Probabilistic Forecasting and {B}ayesian Data Assimilation}.
\newblock Cambridge University Press, 2015.

\bibitem{TdWMR17}
Amirhossein Taghvaei, Jana {de Wiljes}, Prashant~G. Mehta, and Sebastian Reich.
\newblock {{Kalman}} filter and its modern extensions for the continuous-time
  nonlinear filtering problem.
\newblock {\em ASME. J. Dyn. Sys., Meas., Control.}, 140:030904--030904--11,
  2017.

\bibitem{Handel07}
Ramon {van Handel}.
\newblock Stochastic calculus, filtering, and stochastic control.
\newblock Lecture Notes, Caltech, 2007.

\bibitem{villani2021topics}
C{\'e}dric Villani.
\newblock {\em Topics in Optimal Transportation}, volume~58.
\newblock American Mathematical Soc., 2021.

\bibitem{meyn13}
Tao Yang, Prashant~G. Mehta, and Sean~P. Meyn.
\newblock Feedback particle filter.
\newblock {\em IEEE Trans. Automat. Control}, 58(10):2465--2480, 2013.

\bibitem{9811560}
Dongliang Zheng, Jack Ridderhof, Panagiotis Tsiotras, and Ali-akbar
  Agha-mohammadi.
\newblock Belief space planning: {A} covariance steering approach.
\newblock In {\em 2022 International Conference on Robotics and Automation
  (ICRA)}, pages 11051--11057, 2022.

\end{thebibliography}
%
%%%%%%%%%%%%%%%%%%%%%%%%%%%%%%%%

\end{document}